\documentclass[12pt,a4paper]{amsart}
\pdfoutput=1
\usepackage{amsfonts}
\usepackage{amsthm}
\usepackage{amsmath}
\usepackage[utf8]{inputenc}
\usepackage[T1]{fontenc}
\usepackage{lmodern}
\input{glyphtounicode}
\usepackage{indentfirst}
\numberwithin{equation}{section}
\usepackage[margin=2.9cm]{geometry}
\usepackage[hidelinks]{hyperref}

\theoremstyle{plain}
\newtheorem{Th}{Theorem}[section]
\newtheorem{Lemma}[Th]{Lemma}
\newtheorem{Cor}[Th]{Corollary}
\newtheorem{Prop}[Th]{Proposition}

\theoremstyle{definition}
\newtheorem{Def}[Th]{Definition}
\newtheorem{Rem}[Th]{Remark}
\newtheorem{Ex}[Th]{Example}

\newcommand{\E}{\mathrm{E}}
\newcommand{\Var}{\mathrm{Var}}

\begin{document}
	
	\title[Concentration Inequalities for U-Statistics]{Concentration Inequalities for U-Statistics: A Survey with Explicit Constants}
	
	\author[Yannik Pitcan]{Yannik Pitcan}
	
	\address{University of California, Berkeley \\ Department of Statistics \\
		Berkeley CA 94720} 
	
	\email{pitcany@gmail.com}



\begin{abstract} This survey gives a self-contained treatment of concentration inequalities for U-statistics. Using Hoeffding's blocking argument, which reduces a U-statistic to an average over sums of independent random variables, we derive Hoeffding-, Bennett-, and Bernstein-type tail bounds with explicit constants, together with their extensions to unbounded sub-Gaussian and sub-exponential kernels and to two-sample and incomplete U-statistics. The reduction is stated once, as a convex-domination lemma, from which all the tail bounds follow as corollaries. While these results are classical -- the blocking argument goes back to Hoeffding (1963) and Bernstein-type bounds appear in Arcones (1995) -- complete elementary derivations with explicit constants are scattered or omitted in the literature, and collecting them is the purpose of this survey. We close with an overview of sharper bounds available under degeneracy assumptions and of robust median-of-means alternatives for heavy-tailed kernels.
\end{abstract}
	
\maketitle
	
\section{Introduction} U-statistics appear in statistics when producing minimum variance unbiased estimators. Naturally, they appear in clustering, ranking, and learning on graphs. In other words, any time we deal with exchangeable sequences, U-statistics are likely to come up. First, we begin with a definition of a U-statistic following Joly and Lugosi (2016). Hoeffding (1948) was the first to introduce U-statistics in his paper ``A class of statistics with asymptotically normal distribution.''
	
\begin{Def} \label{U-Stat-defined} Let $X$ be a random variable taking values in a measurable space $\mathcal{X}$ and let $h: \mathcal{X}^m\to \mathbb{R}$ be a measurable function of $m\geq 2$ variables. Let $P$ be the probability measure of $X$ and suppose we can access $n\geq m$ i.i.d. random variables $X_1,\ldots,X_n\sim X$. Then the U-statistic of order $m$ and kernel $h$ based on $\{X_i\}$ is defined as
		\[
		U_n(h)=\frac{(n-m)!}{n!} \sum_{(i_1,\ldots,i_m)\in I_n^m} h(X_{i_1},\ldots,X_{i_m}),
		\]	
		where
		\[	
		I_n^m = \{(i_1,\ldots,i_m): 1\leq i_j\leq n,\, i_j\neq i_k \textrm{ if } j\neq k\}
		\] is the set of all $m$-tuples of different integers between $1$ and $n$.
		
		Alternatively, if we restrict $I_n^m$ to be the set of increasing $m$-tuples of integers between $1$ and $n$, i.e. $i_j < i_k$ for $j<k$, then we have 
		\[
		U_n(h)=\frac{1}{\binom{n}{m}} \sum_{(i_1,\ldots,i_m)\in I_n^m} h(X_{i_1},\ldots,X_{i_m}).
		\]
	\end{Def}
	
	U-statistics are unbiased estimators of the mean $m_h = \E h(X_1,\ldots,X_m)$; in the nonparametric setting, where the underlying distribution ranges over all distributions for which $m_h$ is defined, they have minimal variance among all unbiased estimators (Halmos, 1946). A few examples of U-statistics include the sample variance $s_n^2$, which is a U-statistic of order $2$ with kernel given by $h(x,y)=(1/2)(x-y)^2$, and the third $k$-statistic, $ \frac{n}{(n-1)(n-2)} \sum_{i=1}^n (X_i - \bar{X}_n)^3$, which estimates the third cumulant.
	
	It is of interest to understand the concentration of a U-statistic around its expected value. The remainder of the paper is organized as follows. Section 2 states the classical moment generating function tools. Section \ref{sec:main} isolates Hoeffding's blocking argument as a single convex-domination lemma (Lemma \ref{convex-domination}) and derives Hoeffding-, Bennett-, and Bernstein-type bounds for U-statistics from it, with explicit constants throughout, together with remarks comparing the resulting exponents with the bounded-differences approach and with the true variance of $U_n$. Section \ref{sec:unbounded} extends the bounds to unbounded sub-Gaussian and sub-exponential kernels, covering in particular the sample variance of Gaussian data. Section \ref{sec:variants} treats two-sample U-statistics -- including the Mann--Whitney statistic and the AUC -- and incomplete U-statistics. Section \ref{sec:survey} surveys sharper bounds under degeneracy, variance-adaptive and empirical Bernstein inequalities, and robust median-of-means estimators for heavy-tailed kernels. Section \ref{sec:simulation} illustrates numerically how far the explicit bounds sit from the truth, and why.
	
	\section{Classical bounds} We start by taking a look at some classical bounds that do not involve special conditions on the kernel beyond boundedness. Write $R = \sup h - \inf h$ for the range of the kernel and $k = \left\lfloor n/m\right\rfloor$. An argument of Hoeffding (1963, Section 5) establishes that
	\[
	P \left\{ \left| U_n(h)-m_h \right| > R\sqrt{\frac{\log(\frac{2}{\delta})}{2\left\lfloor n/m\right\rfloor }} \right\} \leq \delta
	\] for $h$ bounded, for all $\delta\in(0,1)$, and one can also derive

	\[
	P \left\{ \left| U_n(h)-m_h \right| > \sqrt{\frac{2\sigma^2 \log(\frac{2}{\delta})}{\left\lfloor n/m\right\rfloor }} + \frac{2 b \log(\frac{2}{\delta})}{3\left\lfloor n/m\right\rfloor} \right\} \leq \delta,
	\] where $\sigma^2$ is the variance of $h(X_1,\ldots,X_m)$ and $b = \left\Vert h - m_h\right\Vert_\infty$. These bounds follow from applications of two classic inequalities; we prove them in Section \ref{sec:main} (Corollaries \ref{ustat-hoeffding-delta} and \ref{ustat-bernstein}).
	
	For our purposes, we state the aforementioned inequalities in detail. The proofs can be seen in the book \textit{Concentration Inequalities} by Boucheron, Lugosi, and Massart (2013).
	
	\begin{Lemma}[Hoeffding's lemma]
		Let $X\in [a,b]$ with probability $1$ and $\E X=\mu$.
		Then \[ \E [e^{s(X-\mu)}]\leq e^{s^2(b-a)^2/8}. \]
	\end{Lemma}
	
	\begin{Th}[Hoeffding's inequality] \label{Hoeffding-Inequality} Let $X_1,\ldots,X_n$ be independent observations such that $\E (X_i)=\mu$ and $X_i\in [a,b]$. Then, for any $\epsilon>0$,
		\[
		P(\left|\bar{X}_n - \mu\right|\geq \epsilon)\leq 2e^{-2n\epsilon^2/(b-a)^2}.
		\]
		
	\end{Th}
	
	\begin{Cor} \label{Hoeffding-Inequality-Corollary} If $X_1,\ldots,X_n$ have the same properties as previously (independent, bounded in $[a,b]$, common mean $\mu$), then, with probability at least $1-\delta$,
		\[
		\left|\bar{X}_n - \mu\right|\leq \sqrt{\frac{(b-a)^2}{2n} \log\left(\frac{2}{\delta}\right)}
		\]
		
	\end{Cor}
	
	\begin{Th}[Bernstein's inequality] \label{Bernstein-Inequality}
		Let $X_1,\ldots,X_n$ be independent real-valued random variables with zero mean such that $\left|X_i\right|\leq c$ almost surely. Let $\sigma_i^2=\Var(X_i)$ and $\sigma^2=\frac{1}{n}\sum_{i=1}^n \sigma_i^2.$ Then,
		\[
		P\{\bar{X}_n\geq \epsilon\}\leq \exp\left(-\frac{n \epsilon^2}{2(\sigma^2+c\epsilon/3)}\right).
		\]
	\end{Th}
	
	\begin{Lemma}[Bernstein's lemma]
		Let $X$ be a random variable such that $\E X=\mu$, $\Var(X)=\sigma^2$, and $\left|X-\mu\right|\leq c$ almost surely. Then
		\[ 
		\E [e^{s(X-\mu)}]\leq \exp \left( \frac{\sigma^2}{c^2} (e^{sc}-1-sc)\right). 
		\]
	\end{Lemma}
	
	\section{Hoeffding, Bennett, and Bernstein bounds for U-statistics}\label{sec:main}

We cannot immediately apply the inequalities of Section 2 to $U_n(h)$ because the summands $h(X_{i_1},\ldots,X_{i_m})$ share variables and are therefore dependent. Hoeffding's remedy (1963, Section 5) is to rewrite the U-statistic as a symmetric average of block averages of independent terms. Let $k=\left\lfloor n/m\right\rfloor$ and
\[
V(x_1,\ldots,x_n)=\frac{1}{k} \sum_{i=1}^{k} h(x_{(i-1)m+1},\ldots,x_{im}),
\]
the average of the kernel evaluated on $k$ disjoint blocks. Then
\begin{equation}\label{eq:blocking}
U_n(h) = \frac{1}{n!}\sum_{\sigma\in S_n} V(X_{\sigma(1)},\ldots,X_{\sigma(n)}).
\end{equation}
Indeed, for each fixed block position $j\in\{1,\ldots,k\}$ and each ordered $m$-tuple of distinct indices, exactly $(n-m)!$ permutations $\sigma$ place that tuple in block $j$; since $V$ averages over the $k$ blocks, summing over $\sigma\in S_n$ and dividing by $n!$ recovers the factor $(n-m)!/n!$ in Definition \ref{U-Stat-defined}.

For each fixed $\sigma$, the blocks are disjoint, so $V(X_{\sigma(1)},\ldots,X_{\sigma(n)})$ is distributed as an average of $k$ i.i.d.\ copies of $h(X_1,\ldots,X_m)$. The next lemma turns this representation into a transfer principle: for tail bounds proved through the moment generating function, the U-statistic behaves no worse than an i.i.d.\ average with sample size $k$.

\begin{Lemma}[Convex domination]\label{convex-domination}
	Let $Z_1,\ldots,Z_k$ be i.i.d.\ copies of $h(X_1,\ldots,X_m)$ and $\bar{Z}_k = \frac{1}{k}\sum_{i=1}^k Z_i$. Then for every convex $\Phi:\mathbb{R}\to\mathbb{R}$,
	\[
	\E\, \Phi(U_n(h)-m_h) \leq \E\, \Phi(\bar{Z}_k-m_h).
	\]
	In particular, for all $s\in\mathbb{R}$,
	\[
	\E\, e^{s(U_n(h)-m_h)} \leq \left( \E\, e^{(s/k)(Z_1-m_h)} \right)^k.
	\]
\end{Lemma}

\begin{proof}
	Write $V_\sigma = V(X_{\sigma(1)},\ldots,X_{\sigma(n)})$. By \eqref{eq:blocking} and Jensen's inequality applied to the average over $\sigma$,
	\[
	\Phi(U_n(h)-m_h) = \Phi\left( \frac{1}{n!}\sum_{\sigma\in S_n} (V_\sigma - m_h) \right) \leq \frac{1}{n!}\sum_{\sigma\in S_n} \Phi(V_\sigma - m_h).
	\]
	Taking expectations and using that $V_\sigma$ has the same distribution as $\bar{Z}_k$ for every $\sigma$ gives the first claim. The second claim is the case $\Phi(x)=e^{sx}$, together with independence of the $Z_i$:
	$\E e^{s(\bar{Z}_k - m_h)} = \left(\E e^{(s/k)(Z_1-m_h)}\right)^k$.
\end{proof}

\begin{Th}[Hoeffding's inequality for U-statistics]\label{ustat-hoeffding}
	Suppose $h$ takes values in an interval of length $R$. Then for all $t>0$,
	\[
	P(U_n(h)-m_h\geq t)\leq \exp\left(-\frac{2kt^2}{R^2}\right)
	\quad\text{and}\quad
	P(\left|U_n(h)-m_h\right|\geq t)\leq 2\exp\left(-\frac{2kt^2}{R^2}\right),
	\]
	where $k=\left\lfloor n/m\right\rfloor$.
\end{Th}

\begin{proof}
	The variable $(Z_1-m_h)/k$ takes values in an interval of length $R/k$, so Hoeffding's Lemma gives $\E e^{(s/k)(Z_1-m_h)} \leq e^{s^2R^2/(8k^2)}$. By Lemma \ref{convex-domination},
	\[
	\E\, e^{s(U_n(h)-m_h)} \leq e^{s^2R^2/(8k)},
	\]
	and Markov's inequality yields $P(U_n(h)-m_h\geq t)\leq \exp\left(s^2R^2/(8k) - st\right)$ for every $s>0$. The right-hand side is minimized at $s = 4kt/R^2$, giving $\exp(-2kt^2/R^2)$. Applying the same bound to the kernel $-h$ and a union bound gives the two-sided statement.
\end{proof}

\begin{Cor}\label{ustat-hoeffding-delta}
	Under the assumptions of Theorem \ref{ustat-hoeffding}, for every $\delta\in(0,1)$, with probability at least $1-\delta$,
	\[
	\left|U_n(h)-m_h\right| \leq R\sqrt{\frac{\log(\frac{2}{\delta})}{2\left\lfloor n/m\right\rfloor}}.
	\]
\end{Cor}

\begin{Rem}
	Since $R\leq 2\left\Vert h\right\Vert_\infty$, these bounds also hold with $2\left\Vert h\right\Vert_\infty$ in place of $R$; stating them in terms of the range is sharper (for instance, for a nonnegative kernel $R\leq\left\Vert h\right\Vert_\infty$).
\end{Rem}

We can use the same approach to develop bounds involving the variance of $h(X_1,\ldots,X_m)$.

\begin{Th}[Bennett's inequality for U-statistics]\label{ustat-bennett}
	Let $\sigma^2=\Var(h(X_1,\ldots,X_m))$ and $b=\left\Vert h-m_h\right\Vert_\infty<\infty$. Then for all $t>0$,
	\[
	P(U_n(h)-m_h\geq t)\leq \exp\left(-\frac{k\sigma^2}{b^2}\,\varphi\!\left(\frac{bt}{\sigma^2}\right)\right),
	\qquad \varphi(u)=(1+u)\log(1+u)-u.
	\]
\end{Th}

\begin{proof}
	Apply Bernstein's Lemma to $(Z_1-m_h)/k$, which has variance $\sigma^2/k^2$ and is bounded in absolute value by $c=b/k$:
	\[
	\E\, e^{(s/k)(Z_1-m_h)} \leq \exp\left\{\frac{\sigma^2}{b^2}\left(e^{\frac{sb}{k}}-1-\frac{sb}{k}\right)\right\}.
	\]
	By Lemma \ref{convex-domination} and Markov's inequality, for every $s>0$,
	\[
	P(U_n(h)-m_h\geq t)\leq \exp\left\{\frac{k\sigma^{2}}{b^{2}}\left(e^{\frac{sb}{k}}-1-\frac{sb}{k}\right) - st\right\}.
	\]
	The exponent is minimized at $s=\frac{k}{b}\log(1+\frac{tb}{\sigma^2})$; substituting this choice and writing $u = tb/\sigma^2$, the exponent becomes
	\[
	\frac{k\sigma^2}{b^2}\left(u - \log(1+u) - u\log(1+u)\right) = -\frac{k\sigma^2}{b^2}\,\varphi(u). \qedhere
	\]
\end{proof}

\begin{Cor}[Bernstein's inequality for U-statistics]\label{ustat-bernstein}
	Under the assumptions of Theorem \ref{ustat-bennett}, for all $t>0$,
	\[
	P(\left|U_n(h)-m_h\right|\geq t)\leq 2\exp\left(-\frac{kt^2}{2(\sigma^2+bt/3)}\right),
	\]
	and consequently, for every $\delta\in(0,1)$, with probability at least $1-\delta$,
	\[
	\left|U_n(h)-m_h\right| \leq \sqrt{\frac{2\sigma^2\log(\frac{2}{\delta})}{k}} + \frac{2b\log(\frac{2}{\delta})}{3k}.
	\]
\end{Cor}

\begin{proof}
	The elementary inequality $\varphi(u)\geq \frac{u^2}{2(1+u/3)}$ for $u\geq 0$ turns Bennett's bound into the one-sided Bernstein bound; the two-sided form follows as before by applying the result to $-h$. For the second claim, let $L=\log(\frac{2}{\delta})$. The bound $2\exp(-kt^2/(2(\sigma^2+bt/3)))\leq\delta$ holds if and only if $kt^2 - \frac{2bL}{3}t - 2\sigma^2 L\geq 0$, i.e.\ for
	\[
	t \geq \frac{1}{k}\left(\frac{bL}{3} + \sqrt{\frac{b^2L^2}{9} + 2k\sigma^2 L}\right),
	\]
	and by $\sqrt{x+y}\leq\sqrt{x}+\sqrt{y}$ the right-hand side is at most $\sqrt{2\sigma^2 L/k} + \frac{2bL}{3k}$.
\end{proof}

\begin{Rem}[Comparison with bounded differences]\label{rem-mcdiarmid}
	One might instead apply McDiarmid's bounded-differences inequality directly to $U_n(h)$ as a function of $(X_1,\ldots,X_n)$: changing one coordinate changes $U_n$ by at most $mR/n$, which yields
	$P(\left|U_n(h)-m_h\right|\geq t)\leq 2\exp\left(-2nt^2/(m^2R^2)\right)$.
	The exponent in Theorem \ref{ustat-hoeffding} is $2\left\lfloor n/m\right\rfloor t^2/R^2 \approx 2nt^2/(mR^2)$, larger by a factor of $m$. Thus the blocking argument, although it discards all comparisons across blocks, strictly dominates the generic bounded-differences approach for U-statistics.
\end{Rem}

\begin{Rem}[Sharpness of the variance proxy]\label{rem-variance}
	The effective sample size in these bounds is $k=\left\lfloor n/m\right\rfloor$: the Gaussian-regime term in Corollary \ref{ustat-bernstein} scales as $\sqrt{m\sigma^2/n}$. The actual variance of the U-statistic is $\Var(U_n(h)) = \frac{m^2}{n}\zeta_1 + O(n^{-2})$, where $\zeta_1 = \Var\left(\E[h(X_1,\ldots,X_m)\mid X_1]\right)$. Hoeffding's decomposition gives $\sigma^2\geq m\zeta_1$, with equality precisely when all higher-order projections vanish, so the variance proxy $m\sigma^2/n$ used above is conservative by the factor $\sigma^2/(m\zeta_1)\geq 1$. Bounds that capture the correct $\zeta_1$-dependence require finer arguments; see Arcones (1995).
\end{Rem}

So far, all these results are general and do not involve any assumptions about the kernel $h$ besides measurability and boundedness.

	\section{Unbounded kernels: sub-Gaussian and sub-exponential bounds}\label{sec:unbounded}

The results of Section \ref{sec:main} require a bounded kernel, and this excludes the motivating examples of the introduction whenever the data themselves are unbounded: the sample variance kernel $h(x,y)=(x-y)^2/2$ and the $k$-statistic kernels are unbounded on unbounded domains. Lemma \ref{convex-domination}, however, transfers \emph{any} bound on the moment generating function of the kernel, so light-tailed kernels can be handled by exactly the same argument. We use the standard parametrizations.

\begin{Def}\label{def:subg}
	A real random variable $Z$ with mean $\mu$ is \emph{sub-Gaussian with variance proxy $\nu^2$} if
	\[
	\E\, e^{\lambda(Z-\mu)}\leq e^{\lambda^2\nu^2/2}\quad\text{for all }\lambda\in\mathbb{R},
	\]
	and \emph{sub-exponential with parameters $(\nu^2,\alpha)$}, $\alpha>0$, if the same bound holds for all $\left|\lambda\right|\leq 1/\alpha$.
\end{Def}

\begin{Th}[Sub-Gaussian kernels]\label{ustat-subgaussian}
	Suppose $Z=h(X_1,\ldots,X_m)$ is sub-Gaussian with variance proxy $\nu^2$. Then, with $k=\left\lfloor n/m\right\rfloor$, for all $t>0$,
	\[
	P(\left|U_n(h)-m_h\right|\geq t)\leq 2\exp\left(-\frac{kt^2}{2\nu^2}\right),
	\]
	and for every $\delta\in(0,1)$, with probability at least $1-\delta$,
	\[
	\left|U_n(h)-m_h\right|\leq \nu\sqrt{\frac{2\log(\frac{2}{\delta})}{k}}.
	\]
\end{Th}

\begin{proof}
	By Lemma \ref{convex-domination} with $\lambda=s/k$,
	\[
	\E\, e^{s(U_n(h)-m_h)}\leq \left(e^{(s/k)^2\nu^2/2}\right)^k = e^{s^2\nu^2/(2k)}
	\quad\text{for all }s\in\mathbb{R},
	\]
	i.e.\ $U_n(h)-m_h$ is itself sub-Gaussian with variance proxy $\nu^2/k$. Markov's inequality and optimizing $s=kt/\nu^2$ give the one-sided tail $\exp(-kt^2/(2\nu^2))$; the two-sided bound and the $\delta$-form follow as before.
\end{proof}

\begin{Rem}
	Theorem \ref{ustat-subgaussian} contains Theorem \ref{ustat-hoeffding}: a kernel with range $R$ is sub-Gaussian with $\nu=R/2$ by Hoeffding's Lemma, and $kt^2/(2\nu^2)=2kt^2/R^2$.
\end{Rem}

\begin{Th}[Sub-exponential kernels]\label{ustat-subexponential}
	Suppose $Z=h(X_1,\ldots,X_m)$ is sub-exponential with parameters $(\nu^2,\alpha)$. Then, with $k=\left\lfloor n/m\right\rfloor$, for all $t>0$,
	\[
	P(\left|U_n(h)-m_h\right|\geq t)\leq 2\exp\left(-\frac{k}{2}\min\left\{\frac{t^2}{\nu^2},\,\frac{t}{\alpha}\right\}\right),
	\]
	and for every $\delta\in(0,1)$, with probability at least $1-\delta$,
	\[
	\left|U_n(h)-m_h\right|\leq \nu\sqrt{\frac{2\log(\frac{2}{\delta})}{k}} + \frac{2\alpha\log(\frac{2}{\delta})}{k}.
	\]
\end{Th}

\begin{proof}
	By Lemma \ref{convex-domination} and Markov's inequality, for $0<s\leq k/\alpha$,
	\[
	P(U_n(h)-m_h\geq t)\leq \exp\left(\frac{s^2\nu^2}{2k}-st\right).
	\]
	The unconstrained minimizer is $s^\ast=kt/\nu^2$. If $t\leq\nu^2/\alpha$, then $s^\ast\leq k/\alpha$ is feasible and yields $\exp(-kt^2/(2\nu^2))$. If $t>\nu^2/\alpha$, take $s=k/\alpha$; the exponent is
	\[
	\frac{k\nu^2}{2\alpha^2}-\frac{kt}{\alpha} \leq \frac{kt}{2\alpha}-\frac{kt}{\alpha} = -\frac{kt}{2\alpha},
	\]
	using $\nu^2/\alpha\leq t$. Combining the two regimes gives the one-sided bound $\exp(-\frac{k}{2}\min\{t^2/\nu^2,t/\alpha\})$; the definition is symmetric in $\lambda$, so the same bound holds for $-h$, and a union bound gives the two-sided form. For the $\delta$-form, let $L=\log(\frac{2}{\delta})$ and $t=\nu\sqrt{2L/k}+2\alpha L/k$: then $t\geq \nu\sqrt{2L/k}$ gives $kt^2/\nu^2\geq 2L$ and $t\geq 2\alpha L/k$ gives $kt/\alpha\geq 2L$, so the tail probability is at most $2e^{-L}=\delta$.
\end{proof}

\begin{Ex}[Sample variance of Gaussian data]\label{ex:sample-variance}
	Let $X_1,\ldots,X_n\sim N(\mu,\tau^2)$ and recall that the unbiased sample variance $s_n^2$ is the U-statistic of order $m=2$ with kernel $h(x,y)=(x-y)^2/2$. Here $X_1-X_2\sim N(0,2\tau^2)$, so $Z=h(X_1,X_2)=\tau^2 W$ with $W\sim\chi^2_1$, and the standard $\chi^2$ moment generating function bound $\E\, e^{\lambda(W-1)}\leq e^{2\lambda^2}$ for $\left|\lambda\right|<1/4$ (Wainwright, 2019, Example 2.11) shows that $Z$ is sub-exponential with $(\nu^2,\alpha)=(4\tau^4,4\tau^2)$. Theorem \ref{ustat-subexponential} with $k=\left\lfloor n/2\right\rfloor$ then gives, with probability at least $1-\delta$,
	\[
	\left|s_n^2-\tau^2\right| \leq 2\tau^2\sqrt{\frac{2\log(\frac{2}{\delta})}{\left\lfloor n/2\right\rfloor}} + \frac{8\tau^2\log(\frac{2}{\delta})}{\left\lfloor n/2\right\rfloor}.
	\]
\end{Ex}

\begin{Rem}[Heavier tails]\label{rem-subweibull}
	Theorem \ref{ustat-subexponential} does not cover every polynomial kernel of light-tailed data. For instance, the third $k$-statistic has a kernel containing cubic terms, and the cube of a Gaussian variable has tails of order $\exp(-ct^{2/3})$ -- it is sub-Weibull of order $2/3$, not sub-exponential, and its moment generating function is infinite for every $\lambda\neq 0$, so no bound of the form required by Lemma \ref{convex-domination} applies directly. Nevertheless, Lemma \ref{convex-domination} is not limited to exponential moments: since $x\mapsto\psi(\left|x\right|/c)$ is convex for any convex increasing $\psi$, the lemma yields $\left\Vert U_n(h)-m_h\right\Vert_\psi\leq \left\Vert \bar{Z}_k-m_h\right\Vert_\psi$ for every Orlicz norm. Consequently, generalized Bernstein--Orlicz bounds for averages of sub-Weibull variables (Kuchibhotla and Chakrabortty, 2022) carry over to U-statistics with $n$ replaced by $k=\left\lfloor n/m\right\rfloor$, exactly as in the sub-Gaussian and sub-exponential cases.
\end{Rem}

	\section{Two-sample and incomplete U-statistics}\label{sec:variants}

The transfer principle of Lemma \ref{convex-domination} is not tied to the one-sample setting. In this section we record two variants that arise constantly in applications: two-sample U-statistics, which include the Mann--Whitney statistic and the area under the ROC curve (AUC), and incomplete U-statistics, which are what practitioners actually compute when $\binom{n}{m}$ is prohibitively large.

\subsection{Two-sample U-statistics}

\begin{Def}\label{def:two-sample}
	Let $X_1,\ldots,X_{n_1}$ be i.i.d.\ with law $P$ on $\mathcal{X}$ and, independently, $Y_1,\ldots,Y_{n_2}$ i.i.d.\ with law $Q$ on $\mathcal{Y}$. Given a measurable kernel $h:\mathcal{X}^{m_1}\times\mathcal{Y}^{m_2}\to\mathbb{R}$, the two-sample U-statistic is
	\[
	U_{n_1,n_2}(h) = \frac{1}{\binom{n_1}{m_1}\binom{n_2}{m_2}} \sum_{(i_1,\ldots,i_{m_1})} \sum_{(j_1,\ldots,j_{m_2})} h(X_{i_1},\ldots,X_{i_{m_1}};Y_{j_1},\ldots,Y_{j_{m_2}}),
	\]
	the sums extending over increasing tuples, with mean $m_h = \E\, h(X_1,\ldots,X_{m_1};Y_1,\ldots,Y_{m_2})$.
\end{Def}

\begin{Th}[Convex domination, two-sample case]\label{two-sample-domination}
	Let $k = \min\left(\left\lfloor n_1/m_1\right\rfloor, \left\lfloor n_2/m_2\right\rfloor\right)$ and let $\bar{Z}_k$ be an average of $k$ i.i.d.\ copies of $h(X_1,\ldots,X_{m_1};Y_1,\ldots,Y_{m_2})$. Then for every convex $\Phi:\mathbb{R}\to\mathbb{R}$,
	\[
	\E\, \Phi(U_{n_1,n_2}(h)-m_h) \leq \E\, \Phi(\bar{Z}_k-m_h).
	\]
	Consequently, Theorems \ref{ustat-hoeffding}, \ref{ustat-bennett}, \ref{ustat-subgaussian}, and \ref{ustat-subexponential}, and Corollaries \ref{ustat-hoeffding-delta} and \ref{ustat-bernstein}, all hold for $U_{n_1,n_2}(h)$ with $\left\lfloor n/m\right\rfloor$ replaced by $k$.
\end{Th}

\begin{proof}
	Define the block average
	\[
	V(x_1,\ldots,x_{n_1};y_1,\ldots,y_{n_2}) = \frac{1}{k}\sum_{i=1}^{k} h\!\left(x_{(i-1)m_1+1},\ldots,x_{im_1};\, y_{(i-1)m_2+1},\ldots,y_{im_2}\right),
	\]
	which uses $k$ blocks that are disjoint within each sample (they fit because $km_1\leq n_1$ and $km_2\leq n_2$). For each fixed block position and each pair of ordered tuples, exactly $(n_1-m_1)!\,(n_2-m_2)!$ permutation pairs $(\sigma,\tau)\in S_{n_1}\times S_{n_2}$ place that pair of tuples in the given block, so, as in \eqref{eq:blocking},
	\[
	U_{n_1,n_2}(h) = \frac{1}{n_1!\,n_2!} \sum_{\sigma\in S_{n_1}}\sum_{\tau\in S_{n_2}} V(X_{\sigma(1)},\ldots,X_{\sigma(n_1)};Y_{\tau(1)},\ldots,Y_{\tau(n_2)}).
	\]
	For each fixed $(\sigma,\tau)$ the block evaluations are independent, since they involve disjoint sets of $X$'s and disjoint sets of $Y$'s, so each $V$ term is distributed as $\bar{Z}_k$. Jensen's inequality applied to the average over $(\sigma,\tau)$ concludes, exactly as in the proof of Lemma \ref{convex-domination}. The transferred theorems follow because their proofs used only the moment generating function bound supplied by convex domination.
\end{proof}

\begin{Cor}[Mann--Whitney statistic / AUC]\label{cor:auc}
	Suppose the laws of $X$ and $Y$ are continuous, let $\theta = P(X<Y)$, and let
	\[
	\widehat{\theta}_{n_1,n_2} = \frac{1}{n_1 n_2}\sum_{i=1}^{n_1}\sum_{j=1}^{n_2} \mathbf{1}\{X_i<Y_j\}
	\]
	be the Mann--Whitney statistic (equivalently, the empirical AUC when the $X_i$ are scores of negative examples and the $Y_j$ scores of positive examples). Then for all $t>0$,
	\[
	P\left(\left|\widehat{\theta}_{n_1,n_2}-\theta\right|\geq t\right) \leq 2\exp\left(-2\min(n_1,n_2)\,t^2\right),
	\]
	and with probability at least $1-\delta$,
	\[
	\left|\widehat{\theta}_{n_1,n_2}-\theta\right| \leq \sqrt{\frac{\log(\frac{2}{\delta})}{2\min(n_1,n_2)}}.
	\]
\end{Cor}

\begin{proof}
	This is Theorem \ref{two-sample-domination} combined with Theorem \ref{ustat-hoeffding} for the kernel $h(x;y)=\mathbf{1}\{x<y\}$, which has $m_1=m_2=1$, range $R=1$, and $k=\min(n_1,n_2)$.
\end{proof}

\begin{Rem}\label{rem:min-vs-harmonic}
	The effective sample size $\min(n_1,n_2)$ is the natural one here: since
	\[
	\frac{n_1 n_2}{n_1+n_2} \leq \min(n_1,n_2) \leq \frac{2\,n_1 n_2}{n_1+n_2},
	\]
	it agrees with the harmonic effective sample size $n_1n_2/(n_1+n_2)$, which sets the scale of the variance of $\widehat{\theta}_{n_1,n_2}$, up to a factor of $2$. In the balanced case $n_1=n_2$ nothing is lost. For the use of such bounds in bipartite ranking, where the AUC is the target of empirical risk minimization, see Cl\'emen\c{c}on, Lugosi, and Vayatis (2008).
\end{Rem}

\subsection{Incomplete U-statistics}

Computing $U_n(h)$ requires $\binom{n}{m}$ kernel evaluations, which is often infeasible. A standard remedy is to average over a random subsample of tuples.

\begin{Def}\label{def:incomplete}
	Let $I_n^{m,<}$ denote the set of increasing $m$-tuples from $\{1,\ldots,n\}$. Draw $I_1,\ldots,I_B$ uniformly at random from $I_n^{m,<}$, with replacement and independently of the data. The incomplete U-statistic with $B$ terms is
	\[
	\widetilde{U}_B(h) = \frac{1}{B}\sum_{b=1}^{B} h(X_{I_b}),
	\]
	where $h(X_I)$ denotes the kernel evaluated on the coordinates indexed by $I$.
\end{Def}

Note that $\E[\widetilde{U}_B(h)\mid X_1,\ldots,X_n] = U_n(h)$, so $\widetilde{U}_B$ is unbiased for $m_h$.

\begin{Prop}\label{prop:incomplete}
	Suppose $h$ takes values in an interval of length $R$, and let $k=\left\lfloor n/m\right\rfloor$. For every $\delta\in(0,1)$, with probability at least $1-\delta$,
	\[
	\left|\widetilde{U}_B(h)-m_h\right| \leq R\sqrt{\frac{\log(\frac{4}{\delta})}{2k}} + R\sqrt{\frac{\log(\frac{4}{\delta})}{2B}}.
	\]
\end{Prop}

\begin{proof}
	Decompose $\widetilde{U}_B-m_h = (\widetilde{U}_B-U_n) + (U_n-m_h)$. Conditionally on the data, the variables $h(X_{I_1}),\ldots,h(X_{I_B})$ are i.i.d.\ with mean $U_n(h)$ and values in an interval of length $R$, so Hoeffding's inequality (Theorem \ref{Hoeffding-Inequality}) applied conditionally gives
	\[
	P\left(\left|\widetilde{U}_B-U_n\right|\geq t \,\middle|\, X_1,\ldots,X_n\right) \leq 2e^{-2Bt^2/R^2},
	\]
	and the same bound holds unconditionally by taking expectations. Theorem \ref{ustat-hoeffding} bounds $\left|U_n-m_h\right|$. Choosing $t$ so that each event has probability at most $\delta/2$ and applying the triangle inequality and a union bound gives the claim.
\end{proof}

\begin{Rem}\label{rem:incomplete-B}
	Proposition \ref{prop:incomplete} makes precise the folklore that incompleteness is cheap: taking $B = k$ terms -- a vanishing fraction $k/\binom{n}{m}$ of all tuples -- already matches the rate of the complete U-statistic up to a small constant factor, and $B\gg k$ makes the Monte Carlo term negligible. The same conditional argument yields a Bernstein version, with the conditional variance $\frac{1}{|I_n^{m,<}|}\sum_{I}(h(X_I)-U_n)^2$ in place of $R^2$. Sharper results, including Berry--Esseen bounds and bootstrap approximations for randomized incomplete U-statistics in high dimensions, are given by Chen and Kato (2019); for classical variance computations under deterministic designs, see Blom (1976) and Lee (1990).
\end{Rem}

	\section{Discussion of recent work}\label{sec:survey}

The bounds of Sections \ref{sec:main}--\ref{sec:variants} depend on the kernel only through global quantities -- its range, variance, or moment generating function -- and on the sample only through $k=\left\lfloor n/m\right\rfloor$. Under additional structure, considerably more is known. We organize this section around three themes -- degeneracy, variance adaptivity, and robustness to heavy tails -- closing with a related sampling design.

\subsection{Degenerate kernels}
Arcones and Gin\'e (1993) required that the kernel was symmetric (invariant under permutation of the $X_i$'s) as well as bounded. In the presence of this and a degeneracy assumption, they were able to establish tighter bounds.

\begin{Def}
	A symmetric kernel $h$ is $P$-degenerate of order $q-1,\, 1<q\leq m$, if for all $x_1, \ldots, x_{q-1}\in \mathcal{X}$,

	\[
	\int h(x_1,\ldots,x_m)dP^{m-q+1}(x_q,\ldots,x_m) = \int h(x_1,\ldots,x_m)dP^{m}(x_1,\ldots,x_m)
	\]
	and
	\[
	(x_1,\ldots,x_q)\mapsto \int h(x_1,\ldots,x_m)dP^{m-q}(x_{q+1},\ldots,x_m)
	\]
	is not a constant function (Joly and Lugosi, 2016).
\end{Def}

When the kernel is $P$-degenerate of order $m-1$ and $m_h=\E[h(X_1,\ldots,X_m)]=0$, then $h$ is considered to be $P$-canonical.

The method used by Arcones and Gin\'e to improve concentration guarantees in this scenario relied on decoupling arguments, where they replaced quadratic forms of random variables by bilinear forms. For example, a decoupling argument would be something along the lines of the following, taken from the monograph of Vershynin (2018):

\begin{Ex}
	Let $A$ be an $n\times n$ matrix with zero diagonal and $X=(X_1,\ldots,X_n)$ is a random vector with independent mean 0 coefficients. If $F$ is convex, then
	\[
	\E F(\langle AX,X\rangle)\leq \E F(4\langle AX,X'\rangle)
	\]
	where $X'$ is an independent copy of $X$.
\end{Ex}

Arcones and Gin\'e showed that if $h$ is canonical of $m$ variables, then there exist positive constants $c_1$ and $c_2$ depending only on $m$ such that for all $\delta\in (0,1)$,

\[
P\left\{ \left|U_{n}(h)-m_{h}\right|\geq c_{1}\left\Vert h\right\Vert _{\infty}\left(\frac{\log(\frac{c_{2}}{\delta})}{n}\right)^{m/2}\right\} \leq\delta,
\]
and
\[
P\left\{ \left|U_{n}(h)-m_{h}\right|\geq\left(\frac{\sigma^{2}\log(\frac{c_{1}}{\delta})}{c_{2}n}\right)^{m/2}\lor\frac{\left\Vert h\right\Vert _{\infty}}{\sqrt{n}}\left(\frac{\log(\frac{c_{1}}{\delta})}{c_{2}}\right)^{(m+1)/2}\right\} \leq\delta.
\]

The first of the two proposed inequalities by Arcones and Gin\'e provides a rate of convergence of $O_p(n^{-m/2})$ -- already $O_p(n^{-1})$ for $m=2$ -- which is much faster than the $O_p(n^{-1/2})$ rate implied by the classical Hoeffding bound.

For canonical kernels of order two, Gin\'e, Lata\l{}a, and Zinn (2000) obtained exponential and moment inequalities for unbounded kernels that correctly interpolate between Gaussian, exponential, and heavier tail regimes; Adamczak (2006) extended the moment inequalities to canonical U-statistics of arbitrary order, and Houdr\'e and Reynaud-Bouret (2003) gave versions with explicit constants for order two. All of these results require light tails, however, and fail if the kernel is heavy-tailed.

\subsection{Variance-adaptive bounds for non-degenerate kernels}
Remark \ref{rem-variance} noted that the blocking bounds of Section \ref{sec:main} carry the variance proxy $m\sigma^2/n$, whereas the true variance of a non-degenerate U-statistic is $\frac{m^2}{n}\zeta_1+O(n^{-2})$ with $\zeta_1 = \Var(\E[h(X_1,\ldots,X_m)\mid X_1])$. Arcones (1995) established a Bernstein-type inequality whose Gaussian regime is governed by $\zeta_1$ rather than $\sigma^2$, with constants depending on $m$; this is the natural sharpening of Corollary \ref{ustat-bernstein}. Maurer (2019) reached variance-sensitive Bernstein-type bounds by a different route, quantifying the interaction between the arguments of a general function of independent variables, with non-degenerate U-statistics as a leading example. On the practical side, Peel, Anthoine, and Ralaivola (2010) proved empirical Bernstein inequalities for U-statistics, in which the unknown variance is replaced by an estimate computed from the same sample, making the bounds directly usable at run time.

\subsection{Heavy tails: median-of-means and robust estimators}
For heavy-tailed kernels a different estimator is needed. Nemirovsky and Yudin (1983) introduced the robust median-of-means estimator of the mean of a random variable $X$. What they did was take a confidence level $\delta\in (0,1)$ and divide the data into $V\approx \log \delta^{-1}$ blocks. Then they computed the empirical mean for each block and took the median of these empirical means, call this $\bar{\mu}$. One can show that

\[
\left| \bar{\mu}-\E X \right|\leq c\sqrt{\Var(X)}\sqrt{\frac{\log(1/\delta)}{n}}
\]
with probability at least $1-\delta$ for an absolute constant $c$, requiring only that $\Var(X)$ be finite. This is an improvement in that these guarantees would still hold even if $X$ were heavy-tailed.

Joly and Lugosi (2016) use a median-of-means style estimator as well, except instead of taking the mean on each block, they calculated the empirical means over $m$-tuples of different blocks to compute (decoupled) U-statistics. Then they computed the medians over all these values. Notably, when the kernel was degenerate of order $m-1$, the rate of convergence of their estimator was $(\log (1/\delta)/n)^{m/2}$, which is similar to that of the U-statistics except without any assumption of boundedness. However, the tradeoff is that degeneracy is required. Minsker and Wei (2020) extended this program to operator-valued U-statistics, constructing robust modifications with exponential deviation guarantees in operator norm and applying them to covariance estimation under the existence of fourth moments only.

\subsection{Sampling without replacement}
Finally, a companion literature studies U-statistics computed on samples drawn without replacement from a finite population. Ai, Ku\v{z}elka, and Wang (2022) prove Hoeffding- and Bernstein-type inequalities in this setting, paralleling the results of Section \ref{sec:main}; the comparison between sampling with and without replacement goes back to Hoeffding (1963), whose convex-ordering argument for that problem is a close relative of Lemma \ref{convex-domination}.

	\section{A numerical illustration}\label{sec:simulation}

We close the technical development with a small simulation that makes the remarks of Section \ref{sec:main} quantitative. Take $X\sim \mathrm{Uniform}[0,1]$ and consider two kernels of order $m=2$:
\begin{itemize}
	\item[(a)] the non-degenerate sample variance kernel $h_a(x,y)=(x-y)^2/2$, for which $m_h=1/12$, $R=1/2$, and $b=\left\Vert h_a-m_h\right\Vert_\infty = 5/12$. Since $D=X-Y$ has the triangular density $1-\left|d\right|$ on $[-1,1]$, $\E D^4 = 1/15$, whence $\sigma^2 = \tfrac14\Var(D^2) = 7/720$; and $\E[h_a\mid X_1=x] = \tfrac12\left((x-\tfrac12)^2+\tfrac1{12}\right)$ gives $\zeta_1 = \tfrac14\Var\left((X-\tfrac12)^2\right) = 1/720$. The variance-proxy gap of Remark \ref{rem-variance} is therefore $\sigma^2/(m\zeta_1) = 7/2$.
	\item[(b)] the $P$-canonical kernel $h_b(x,y)=(x-\tfrac12)(y-\tfrac12)$, for which $m_h=0$, $b=1/4$, $\sigma^2 = 1/144$, and $\zeta_1=0$.
\end{itemize}
For each $n\in\{50,\,200,\,1000,\,5000\}$ we draw $2\times 10^5$ independent samples of size $n$, compute $U_n$ for both kernels, and record the empirical $95\%$ quantile of $\left|U_n-m_h\right|$. Table \ref{tab:simulation} compares these with the Hoeffding bound (Corollary \ref{ustat-hoeffding-delta}), the Bernstein bound (Corollary \ref{ustat-bernstein}), both at $\delta=0.05$, and, for the non-degenerate kernel, the normal approximation $z_{0.975}\sqrt{m^2\zeta_1/n}$ suggested by the asymptotic variance. The simulation script, with a fixed random seed, is provided as an ancillary file.

\begin{table}[ht]
	\centering
	\begin{tabular}{r cccc c cc}
		\hline
		& \multicolumn{4}{c}{kernel (a), non-degenerate} & & \multicolumn{2}{c}{kernel (b), canonical} \\
		\cline{2-5}\cline{7-8}
		$n$ & empirical & Bernstein & Hoeffding & normal & & empirical & Bernstein \\
		\hline
		$50$   & $0.0212$ & $0.0946$ & $0.1358$ & $0.0207$ & & $0.00480$  & $0.0699$ \\
		$200$  & $0.0104$ & $0.0370$ & $0.0679$ & $0.0103$ & & $0.00118$  & $0.0288$ \\
		$1000$ & $0.0046$ & $0.0140$ & $0.0304$ & $0.0046$ & & $0.00024$  & $0.0114$ \\
		$5000$ & $0.0021$ & $0.0058$ & $0.0136$ & $0.0021$ & & $0.00005$  & $0.0048$ \\
		\hline
	\end{tabular}
	\medskip
	\caption{Empirical $95\%$ quantiles of $\left|U_n-m_h\right|$ over $2\times10^5$ replications, against the bounds of Section \ref{sec:main} at $\delta=0.05$ and the normal approximation $z_{0.975}\sqrt{m^2\zeta_1/n}$.}
	\label{tab:simulation}
\end{table}

Three observations. First, for the non-degenerate kernel the normal approximation built on $\zeta_1$ matches the empirical quantile to within two percent already at $n=50$: the projection variance, not $\sigma^2$, is what governs the actual concentration of $U_n$. Second, the Bernstein bound is valid but conservative by a stable factor -- about $2.8$ at $n=5000$. This factor is fully explained: asymptotically it equals
\[
\sqrt{\frac{\sigma^2}{m\zeta_1}} \times \frac{\sqrt{2\log(2/\delta)}}{z_{0.975}} = \sqrt{3.5}\times 1.39 \approx 2.6,
\]
the product of the variance-proxy gap of Remark \ref{rem-variance} and the usual slack of a sub-Gaussian tail bound against the central limit theorem at moderate $\delta$. The Hoeffding bound pays an additional range-versus-variance premium and sits at a stable factor of about $6.6$. Third, for the canonical kernel the picture changes qualitatively: the empirical quantile decays at the degenerate rate $\Theta(n^{-1})$ while the bounds of Section \ref{sec:main} can only see the $\Theta(n^{-1/2})$ scale, so the ratio of bound to truth grows like $\sqrt{n}$ -- from about $15$ at $n=50$ to about $100$ at $n=5000$. This is precisely the gap closed by the degeneracy-aware inequalities surveyed in Section \ref{sec:survey}.

	\section{Concluding remarks} We have given a self-contained account, with explicit constants, of the standard Hoeffding-, Bennett-, and Bernstein-type concentration bounds for U-statistics. The underlying blocking argument -- ``creating independence where there is not'' by representing the U-statistic as a symmetric average over sums of independent blocks -- is due to Hoeffding (1963, Section 5), and Bernstein-type inequalities for U-statistics appear in Arcones (1995); our aim here has been expository, isolating the reduction as a single convex-domination lemma so that each tail bound for sums of independent variables transfers immediately to U-statistics -- including sub-Gaussian and sub-exponential kernels, the two-sample case, and, through a conditional argument, incomplete U-statistics. We then summarized the tradeoffs among current approaches to moment and tail bounds for U-statistics, with emphasis on the techniques of decoupling and median-of-means estimators.


\begin{thebibliography}{99}

		\bibitem{adamczak}
		Adamczak, R. (2006),
		\textit{Moment inequalities for U-statistics.}
		The Annals of Probability 34 (6), 2288-2314

		\bibitem{ai-kuzelka-wang}
		Ai, J., Ku\v{z}elka, O., and Wang, Y. (2022),
		\textit{Hoeffding and Bernstein inequalities for U-statistics without replacement.}
		Statistics \& Probability Letters 187

		\bibitem{arcones1995}
		Arcones, M. A. (1995),
		\textit{A Bernstein-type inequality for U-statistics and U-processes.}
		Statistics \& Probability Letters 22 (3), 239-247

		\bibitem{arcones-gine}
		Arcones, M. A. and Gin\'e, E. (1993),
		\textit{Limit theorems for U-processes.}
		The Annals of Probability 21 (3), 1494-1542

		\bibitem{blom}
		Blom, G. (1976),
		\textit{Some properties of incomplete U-statistics.}
		Biometrika 63 (3), 573-580

		\bibitem{concentration}
		Boucheron, S., Lugosi, G., and Massart, P. (2013),
		\textit{Concentration Inequalities: A Nonasymptotic Theory of Independence.}
		Oxford University Press.

		\bibitem{chen-kato}
		Chen, X. and Kato, K. (2019),
		\textit{Randomized incomplete U-statistics in high dimensions.}
		The Annals of Statistics 47 (6), 3127-3156

		\bibitem{clemencon}
		Cl\'emen\c{c}on, S., Lugosi, G., and Vayatis, N. (2008),
		\textit{Ranking and empirical minimization of U-statistics.}
		The Annals of Statistics 36 (2), 844-874

		\bibitem{decoupling}
		de la Pe\~na, V. H. and Gin\'e, E. (1999),
		\textit{Decoupling. From dependence to independence. Randomly stopped processes. U-statistics and processes. Martingales and beyond.}
		Probability and its Applications (New York). Springer-Verlag, New York.

		\bibitem{gine-latala-zinn}
		Gin\'e, E., Lata\l{}a, R., and Zinn, J. (2000),
		\textit{Exponential and moment inequalities for U-statistics.}
		In: High Dimensional Probability II, Progress in Probability 47, Birkh\"auser, Boston, 13-38

		\bibitem{halmos}
		Halmos, P. R. (1946),
		\textit{The theory of unbiased estimation.}
		The Annals of Mathematical Statistics 17, 34-43

		\bibitem{hoeffding2}
		Hoeffding, W. (1948),
		\textit{A class of statistics with asymptotically normal distribution.}
		The Annals of Mathematical Statistics 19 (3), 293-325

		\bibitem{hoeffding}
		Hoeffding, W. (1963),
		\textit{Probability inequalities for sums of bounded random variables.}
		Journal of the American Statistical Association 58 (301), 13-30

		\bibitem{houdre-rb}
		Houdr\'e, C. and Reynaud-Bouret, P. (2003),
		\textit{Exponential inequalities, with constants, for U-statistics of order two.}
		In: Stochastic Inequalities and Applications, Progress in Probability 56, Birkh\"auser, Basel, 55-69

		\bibitem{lugosi}
		Joly, E. and Lugosi, G. (2016),
		\textit{Robust estimation of U-statistics.}
		Stochastic Processes and their Applications 126 (12), 3760-3773

		\bibitem{kuchibhotla}
		Kuchibhotla, A. K. and Chakrabortty, A. (2022),
		\textit{Moving beyond sub-Gaussianity in high-dimensional statistics: applications in covariance estimation and linear regression.}
		Information and Inference: A Journal of the IMA 11 (4), 1389-1456

		\bibitem{lee}
		Lee, A. J. (1990),
		\textit{U-Statistics: Theory and Practice.}
		Marcel Dekker, New York.

		\bibitem{maurer}
		Maurer, A. (2019),
		\textit{A Bernstein-type inequality for functions of bounded interaction.}
		Bernoulli 25 (2), 1451-1471

		\bibitem{minsker-wei}
		Minsker, S. and Wei, X. (2020),
		\textit{Robust modifications of U-statistics and applications to covariance estimation problems.}
		Bernoulli 26 (1), 694-727

		\bibitem{nemirovsky-yudin}
		Nemirovsky, A. S. and Yudin, D. B. (1983),
		\textit{Problem Complexity and Method Efficiency in Optimization.}
		Wiley, New York.

		\bibitem{peel}
		Peel, T., Anthoine, S., and Ralaivola, L. (2010),
		\textit{Empirical Bernstein inequalities for U-statistics.}
		In: Advances in Neural Information Processing Systems 23

		\bibitem{vershynin}
		Vershynin, R. (2018),
		\textit{High-Dimensional Probability: An Introduction with Applications in Data Science.}
		Cambridge University Press.

		\bibitem{wainwright}
		Wainwright, M. J. (2019),
		\textit{High-Dimensional Statistics: A Non-Asymptotic Viewpoint.}
		Cambridge University Press.
	\end{thebibliography}
\end{document}